\documentclass[12pt]{article}

\usepackage{amssymb}
\usepackage{newtxtext}
\usepackage{amsmath,amsfonts,amsthm,amssymb}
\usepackage[utf8]{inputenc}
\usepackage[T1]{fontenc}    
\usepackage{url}        
\usepackage{geometry}
\usepackage{booktabs}       
\usepackage{amsfonts}       
\usepackage{nicefrac}       
\usepackage{microtype}      
\usepackage{xcolor}         
\usepackage{subfigure}
\usepackage[inline]{enumitem}
\usepackage{cases}
\usepackage{algorithm}
\usepackage{algpseudocode}
\usepackage{graphicx}
\usepackage{dblfloatfix}
\usepackage{float}
\theoremstyle{plain}
\usepackage{makecell}
\usepackage{bm}
\usepackage{setspace}
\usepackage{thmtools}
\usepackage{thm-restate}
\usepackage{cases}
\usepackage{empheq}
\definecolor{dark-gray}{gray}{0.3}
\definecolor{dkgray}{rgb}{.4,.4,.4}
\definecolor{dkblue}{rgb}{0,0,.5}
\definecolor{medblue}{rgb}{0,0,.75}
\definecolor{rust}{rgb}{0.5,0.1,0.1}
\definecolor{darkblue}{rgb}{0,0.08,0.45}
\usepackage[colorlinks]{hyperref}
\usepackage[backend=biber]{biblatex}
\hypersetup{urlcolor=rust}
\hypersetup{citecolor=darkblue}
\hypersetup{linkcolor=blue}

\newcommand{\Diag}{\mathop{\bf diag}}

\newcommand{\one}{{\mathbf{1}}}

\numberwithin{theorem}{section}

\numberwithin{remark}{section}

\numberwithin{definition}{section}

\numberwithin{assumption}{section}

\declaretheorem[name=Theorem,numberwithin=section]{thm-rest}
\declaretheorem[name=Lemma,numberwithin=section]{lem-rest}
\declaretheorem[name=Corollary,numberwithin=section]{cor-rest}
\declaretheorem[name=Definition,numberwithin=section]{def-rest}
\declaretheorem[name=Proposition,numberwithin=section]{prop-rest}
\declaretheorem[name=Assumption,numberwithin=section]{ass-rest}

\allowdisplaybreaks

\DeclareMathAlphabet{\mathcalorigin}{OMS}{cmsy}{m}{n}

\usepackage[colorinlistoftodos,prependcaption,backgroundcolor=black!5!white,bordercolor=red]{todonotes}
\usepackage[capitalise]{cleveref}
\title{\vspace{-2.0cm} Convex Modeling of Price Cross-Impact over Time}
\author{
Vincent Yinjun-Wang\footnote{Stanford MS\&E. Email: yinjunw@stanford.edu}
\and
Madeleine Udell\footnote{Stanford MS\&E. Email: udell@stanford.edu}
}
\date{}
\begin{document}
\maketitle
\begin{abstract}
Transaction costs can make or break a trading strategy, particularly
in relative-value trading of commodity and macro markets, where
edges are a few basis points.  Price impact is a central component of
transaction cost.  Price impact models usually
include self-impact (a trade in a contract moves that contract's
price) but omit two well-documented effects: cross-impact (a trade in
one contract also moves the prices of related contracts) and
transient impact (price impact decays over time, so an unwind
recovers part of the entry cost).  A model without these effects overprices the impact of relative-value trades, whose correlated legs
are built and unwound over days, and so forgoes potentially
profitable trades.  This paper models both effects with a convex quadratic
cost.  In each period, a positive semidefinite matrix built from
volatility, volume, and correlation forecasts couples trades across
contracts.  A power-law decay kernel then couples trades across
periods.  The resulting cost admits no price manipulation even when
liquidity varies over the planning horizon. The model is demonstrated empirically on calendar spread trading of crude
oil futures around the commodity index roll.
\end{abstract}
\newpage
\section{Introduction}\label{s-intro}
 
Multi-period portfolio construction chooses both what to trade and when to
trade it.  We take two empirical regularities as requirements for a useful transaction cost model.  First, trades in related contracts\footnote{We use the
term ``contract'' to refer to tradable assets and derivatives throughout
this paper.} interact: a relative-value trade is cheaper than the sum of its
independently executed legs.
Second, price impact persists: a later trade in the same direction pays the
residual impact of earlier trades as an additional cost, while an unwind
recovers part of the entry cost.  These effects are usually called
\emph{cross-impact} and \emph{transient impact}. Our model captures both
simultaneously.
 
Modeling these two impacts is especially useful when liquidity varies predictably
over the planning horizon. Published roll calendars, scheduled reports,
settlement windows, and expiry cycles all generate recurring changes in volume
and volatility. Relative-value trading in macro and commodity markets is the core
application. Calendar spread trading of crude oil futures provides a
running example, which we study empirically in
Section~\ref{s-experiments}.
 
\paragraph{The index roll.}
Commodity index exposure (e.g., to the S\&P GSCI) is typically
implemented through long positions in commodity futures. Because futures contracts expire, the index periodically rolls its positions by selling the expiring nearby contract and buying the next deferred contract before physical delivery. This scheduled, passive order flow creates temporary selling pressure in the nearby contract and buying pressure in the deferred contract, narrowing the nearby-minus-deferred calendar spread. Once the roll is complete, the spread tends to recover as the price pressure dissipates. From 2004 to 2011, the WTI nearby-minus-deferred spread fell by roughly 30--40 basis points around the roll window and recovered within about two weeks \cite{irwin2023orderflow}.  The index roll schedule is public and fixed in advance. For example, S\&P GSCI rolls one fifth of its position on each of the fifth through ninth business days of the month \cite{spdji2026gsci}. The approximate size of the flow can also be inferred from public position reports and past roll activity.   
 
This predictable spread compression creates a speculative trade: build a short spread position (short the nearby contract and long the deferred contract) before the roll, then reverse it (long the nearby, short the deferred) into the recovery \cite{mou2011limits}.
Recognizing this price dislocation, however, does not
determine how early to enter, how much to trade each day, or when to
unwind.  Intuitively, greater liquidity means lower transaction
costs, so traders have an incentive to shift their trades toward
liquid periods.  Making these choices by discretion is further
complicated by event-driven liquidity changes: for example, the
EIA's Weekly Petroleum Status Report, released Wednesdays at 10:30
a.m.\ ET, concentrates liquidity on the release day
\cite{eia2026wpsr,bjursell2015inventory}.  Existing multi-period
portfolio construction frameworks \cite{boyd2017multiperiod} model
the price impact cost as separable, contract by contract and period
by period, which overprices exactly relative-value trades:
correlated legs built and unwound over days, where the legs' impacts
partially cancel and the entry's impact is still alive at the exit.

\paragraph{Contributions.}
We propose a convex quadratic model of transient cross-impact cost.
At time period \(t\), a positive semidefinite matrix \(\mathbf{\Lambda}_t\) describes
self- and cross-impact under that period's liquidity forecast. A power-law
decay kernel then couples trades across time periods. We show that the resulting
impact cost is convex. Convexity simplifies
computation, and encodes a no-arbitrage condition. A nonconvex impact
cost can assign negative expected cost to some round trip (a trade that
buys and then sells back to zero inventory), and exploiting such
\emph{price manipulation} is illegal \cite{huberman2004price,gatheral2010noarb}.
Requiring convexity rules it out of the model. We empirically show the overpricing of existing
models, and the usefulness of ours, on public data, trading calendar
spreads of crude oil futures.

\section{Related work}\label{s-related}

\paragraph{Portfolio construction and execution.}
Multi-period portfolio construction models trade expected return against risk and
transaction cost over a planning horizon
\cite{boyd2017multiperiod,boyd2024markowitz}. Their price impact cost is modeled as separable,
contract by contract and period by period.  In
Section~\ref{s-experiments}, we empirically study how such separable
modeling overprices the impact and yields a highly conservative
portfolio.

Limit-order-book and optimal execution models
\cite{obizhaeva2013optimal,bouchaud2018trades,jaber2024optimal} address the
complementary problem of executing each period-level trade. Our model operates at portfolio construction timescale.

\paragraph{Cross impact.}
Price impact across related contracts is often aligned with the factor structure
of contracts' returns \cite{benzaquen2017dissecting,bucci2020coimpact,
mastromatteo2017trading}.  Following this observation, we use return
correlation as a proxy for impact correlation.  Many models assume
cross-impact proportional to return covariance (the risk matrix),
\(\mathbf{\Lambda}_t=\lambda\mathbf{\Sigma}_t\)
\cite{garleanu2013dynamic,garleanu2016dynamic}.  Our construction
(Section~\ref{s-inst}) keeps the shared correlation matrix but allows the
cross-impact and risk matrices different scales.

\paragraph{Transient impact.}
Empirical studies find that price impact decays according to a power law, more slowly than an
exponential \cite{bouchaud2003fluctuations,bouchaud2009markets,
bucci2018slow,bouchaud2018trades}.  Exponential decay is nonetheless common
despite its empirical defects, because the impact state is then Markovian and closed-form decision policy can be derived
\cite{obizhaeva2013optimal,
garleanu2013dynamic}.  Power-law decay is not Markovian, so we optimize the
sequence of all trades jointly under the predicted liquidity over the
planning horizon.

\paragraph{Linear and square-root price impact.} Empirically, the single-contract price impact of
a metaorder of size $Q$ (a succession of trades, all executed in the
same direction and originating from
the same market participant) scales as the square root of \(Q\) \cite{loeb1983trading,almgren2005direct,
toth2011anomalous,
moro2009market}. In this paper, we do not reproduce the square-root law, but model single-contract, single-period impact as linear with power-law
decay. Nonlinear impact with decay generally admits illegal price manipulation \cite{gatheral2010noarb,schneider2019cross}, while linear impact with a positive definite decay kernel rules it out
\cite{alfonsi2012order,alfonsi2016multivariate}.

\paragraph{Convexity.} For linear price impact and resulting quadratic cost, positive
semidefiniteness of the full matrix over contracts and times is a sufficient condition for the absence of negative-cost round trips, i.e., price manipulation is never beneficial
\cite{huberman2004price,gatheral2010noarb}.  A collection of positive
semidefinite single-period matrices does not by itself ensure this property
when liquidity changes over time \cite{fruth2014optimal}.

\section{A multi-period cross-impact model}\label{s-model}
In this section, we build the cost of transient cross-impact \(\phi^{\mathrm{imp}}\) in two steps. First, inspired by risk matrix, we introduce a
cross-impact matrix for each time period, and couple the periods through a temporal decay kernel.

We trade \(n\) contracts at dates
\(a_1<\cdots<a_T\). The trade at date \(a_t\) executes over period \(t\).  The
vector \(\mathbf{u}_t\in\mathbb{R}^n\) is the signed dollar
value traded at date \(t\), and dollar holdings satisfy
\[
  \mathbf{x}_t=\mathbf{x}_{t-1}+\mathbf{u}_t,
  \qquad t=1,\ldots,T.
\]
The dates need not be equally spaced. They can be days, or intraday windows
around scheduled events such as an index roll or a storage report. Write \(\mathbf{u}=\big(\mathbf{u}_1,\dots,\mathbf{u}_T\big)\in\mathbb{R}^{nT}\).

\subsection{Single-period cross-impact}\label{s-inst}
We first model price impact of trading a single contract \(i\) at a single period \(t\), under that period's liquidity
forecast.
Following \cite{boyd2017multiperiod}, let \(v_{t,i}>0\) and \(\sigma_{t,i}>0\) denote forecasted dollar volume and
return volatility for contract \(i\) over period \(t\).  
A trade moves the price more when the contract is volatile, and less when it is
heavily traded.  We model the self-impact of contract \(i\) at date \(t\) as
\begin{equation}\label{e-single-impact}
    \gamma_{t,i} u_{t,i},\qquad\text{where }\gamma_{t,i}=\ell_i\frac{\sigma_{t,i}}{v_{t,i}},
\end{equation}
where \(\ell_i \in \mathbb{R}\) is a learnable parameter.  

Trading one contract also moves the prices of related contracts, so we combine
the single-contract price impact \(\gamma_{t,i}\) with a correlation across contracts.  Define
\(\mathbf{\Gamma}_t=\Diag(\gamma_{t,1},\ldots,\gamma_{t,n})\), and let
\(\mathbf{C}_t\in\mathbb{R}^{n\times n}\) be a positive semidefinite matrix with unit
diagonal. We model the cross-impact over period \(t\) as
\begin{equation*}\label{e-lambda}
  \mathbf{\Lambda}_t
  =\mathbf{\Gamma}_t^{1/2}\mathbf{C}_t\mathbf{\Gamma}_t^{1/2},
\end{equation*}
so that the cost of such cross-impact is \(\mathbf{u}_t^\top \mathbf{\Lambda}_t \mathbf{u}_t\).
Because $\mathbf{\Gamma}_t$ is diagonal and \(\mathbf{C}_t\) has unit diagonal, the construction of \(\mathbf{\Lambda}_t\) adds cross-impact without
disturbing the self-impact \eqref{e-single-impact}:
\(\mathbf{\Lambda}_t\) has diagonal
\(\gamma_{t,1},\ldots,\gamma_{t,n}\). We refer to \(\mathbf{C}_t\) as impact correlation, and use return correlation as its proxy \cite{garleanu2013dynamic,garleanu2016dynamic}. The cross-impact matrix and the risk matrix share the same cross-sectional structure, as the risk matrix can be decomposed as 
\[
\mathbf{\Sigma}_t^{1/2} \mathbf{C}_t \mathbf{\Sigma}_t^{1/2},
\]
where \(\mathbf{\Sigma}_t^{1/2}=\Diag(\sigma_{t,1},\ldots,\sigma_{t,n})\). In practice, \(\mathbf{C}_t\) can be further calibrated using limit order book data \cite{benzaquen2017dissecting,bucci2020coimpact,
mastromatteo2017trading}.

\subsection{Transient impact}\label{s-transient}
Cross-impact does not vanish when the period ends. We now couple
cross-impact over time.  We introduce a temporal decay kernel
\begin{equation*}\label{e-kernel}
  k(h)=(1+h/\tau)^{-\beta},
  \qquad \tau,\beta>0,
\end{equation*}
consistent with the power-law decay found empirically
\cite{bouchaud2003fluctuations,bucci2018slow}.  The kernel \(k(\cdot):[0,\infty)\to[0,1]\) is nonincreasing, with \(k(0)=1\). At date \(a_t\), price impact of a trade at date \(a_s\) has decayed to a fraction \(k(a_t-a_s)\) of its initial size. 
The width \(\tau\) sets the timescale of the decay: larger \(\tau\), slower decay.

To charge each trade against the decayed impact of earlier trades, we model the impact cost of trades \(\mathbf{u}\) as
\begin{equation}\label{e-causal-cost}
  \phi^{\mathrm{imp}}(\mathbf{u})
  =\sum_{t=1}^{T}\sum_{s\leq t}k(a_t-a_s)\mathbf{u}_t^\top\mathbf{\Lambda}_t^{1/2}\mathbf{\Lambda}_s^{1/2}\mathbf{u}_s,
\end{equation}
with \(k(0)=1\).
This cost is a convex quadratic.  To  write the cost more compactly, we collect the kernel values in
\(\mathbf{K}\in\mathbb{R}^{T\times T}\):
\[
K_{ts}
=
\begin{cases}
\displaystyle k\left(a_t-a_s\right),
& s\leq t,\\[6pt]
0, & s>t.
\end{cases}
\]
\(\mathbf{K}\) is lower triangular to enforce causality.
The quadratic form is unchanged if we replace \(\mathbf{K}\) by symmetric matrix
\begin{equation*}\label{e-G}
  \mathbf{G}=\frac{\mathbf{K} + \mathbf{K}^\top}{2}.
\end{equation*}
The impact cost of trades \(\mathbf{u}\) is then
\[
\phi^{\mathrm{imp}}(\mathbf{u})=\mathbf{u}^\top\mathbf{A}\mathbf{u},
\] 
where
\begin{equation}\label{e-sandwich}
  \mathbf{B}
  =\Diag(\mathbf{\Lambda}_1^{1/2},\ldots,
         \mathbf{\Lambda}_T^{1/2}),
  \qquad
  \mathbf{A}
  =\mathbf{B}(\mathbf{G}\otimes\mathbf{I}_n)\mathbf{B}.
\end{equation}
In block form,
\begin{equation}\label{e-A}
  \mathbf{A}_{ts}
  =G_{ts}\mathbf{\Lambda}_t^{1/2}\mathbf{\Lambda}_s^{1/2}.
\end{equation}
And in particular, \(\mathbf{A}_{tt}=\mathbf{\Lambda}_t\), recovering the single-period model of Section~\ref{s-inst}.

Our model uses the symmetric coupling \(\mathbf{\Lambda}_t^{1/2}\mathbf{\Lambda}_s^{1/2}\), which guarantees convexity. In contrast, the true, causal cost with time-varying liquidity, $G_{ts}\mathbf{\Lambda}_s$, can make \(\mathbf{A}\) indefinite and admit negative-cost
round trips, i.e., illegal price manipulation.  Appendix~\ref{a-arith} gives explicit examples.  The symmetric coupling avoids this failure.

\subsection{Convexity and no price manipulation}\label{s-convexity}
With linear price impact, positive semidefiniteness of \(\mathbf{A}\) is
sufficient for the absence of negative-cost round trips: it gives
\[
\phi^{\mathrm{imp}}(\mathbf{u})=\mathbf{u}^\top\mathbf{A}\mathbf{u}\geq0
\quad\text{for every }\mathbf{u},
\]
so the model admits no price manipulation.  

We now show that \(\mathbf{A}\succeq0\) holds for any
\(\mathbf{\Lambda}_t\succeq0\) constructed in Section~\ref{s-inst},
so the cost \eqref{e-causal-cost} is convex no matter how liquidity
varies over the horizon.  The argument has two steps: we write the
power-law kernel as a mixture of exponential kernels with nonnegative
weights, and each exponential kernel is positive semidefinite because
it is the characteristic function of a Cauchy distribution.

Define \(\mathbf{H}\in\mathbb{R}^{T\times T}\) by
\[
H_{ts}=k(|a_t-a_s|),
\qquad t,s=1,\ldots,T.
\]
\emph{Step 1: the kernel is a nonnegative mixture of exponentials.}
The substitution \(w=r(1+h/\tau)\) in
\(\Gamma(\beta)=\int_0^\infty w^{\beta-1}e^{-w}\,dw\) gives
\(\int_0^\infty r^{\beta-1}e^{-r(1+h/\tau)}\,dr=\Gamma(\beta)\,k(h)\),
so that
\[
H_{ts}
=
\frac{1}{\Gamma(\beta)}
\int_0^\infty
r^{\beta-1}e^{-r}\,
e^{-(r/\tau)|a_t-a_s|}\,dr .
\]
\emph{Step 2: each exponential kernel is positive semidefinite.}
For \(c>0\), the Cauchy density with scale \(c\) has characteristic
function \(e^{-c|h|}\),
\[
e^{-c|h|}
=\frac{1}{\pi}\int_{-\infty}^{\infty}
\frac{c}{c^{2}+\omega^{2}}\,e^{\mathrm{i}\omega h}\,d\omega ,
\]
so for every \(\mathbf{z}\in\mathbb{R}^{T}\),
\[
\sum_{t,s}z_t z_s\,e^{-c|a_t-a_s|}
=\frac{1}{\pi}\int_{-\infty}^{\infty}
\frac{c}{c^{2}+\omega^{2}}
\Big|\sum_{t=1}^{T}z_t e^{\mathrm{i}\omega a_t}\Big|^{2}\,d\omega
\;\geq\;0 .
\]
Combining the two steps, by taking \(c=r/\tau\) in Step~2 and
integrating against the nonnegative weight
\(r^{\beta-1}e^{-r}/\Gamma(\beta)\) from Step~1, gives
\(\mathbf{H}\succeq0\).  Therefore
\(\mathbf{G}\) is positive semidefinite,
\[
\mathbf{G}
=
\frac{\mathbf{K}+\mathbf{K}^\top}{2}
=
\frac12(\mathbf{I}_T+\mathbf{H})
\succeq \frac12\mathbf{I}_T ,
\]
and by \eqref{e-sandwich} so is \(\mathbf{A}\),
\[
\mathbf{A}\succeq
\frac12\Diag(\mathbf{\Lambda}_1,\ldots,\mathbf{\Lambda}_T)
\succeq0 .
\]
Every \(\mathbf{u}\) thus costs at least half the sum of its single-period
costs, and no round trip has negative cost.

\section{Numerical study}\label{s-experiments}

We evaluate the model by relative-value trading of WTI calendar
spreads around the S\&P GSCI index roll.  The study runs on public data.  While a better practice of 
speculative traders would certainly require licensed exchange data,
the study shows how our careful model of cross-impact and transient impact
changes the constructed portfolio, compared with existing baselines.

\subsection{WTI index-roll positioning}\label{s-exp-roll}    

\paragraph{Setup.}
In each month we trade \(n=2\) contracts:
the nearby contract the index is selling and the deferred contract it
is buying.  The S\&P GSCI roll window is business days five
through nine. Trading takes place on the business days of the month, and
the position is closed by the business day before the nearby
contract's last trading day.

Let \(\mathbf{r}_t\in\mathbb{R}^2\) denote the forecast returns of
the two contracts over period \(t\).  The trader chooses, over the month, how large a spread position to hold and when to build and
unwind it, by solving
\begin{equation}\label{e-roll}
\begin{array}{ll}
\mbox{maximize} &
\displaystyle
\sum_{t=1}^T
\mathbf{r}_t^\top\mathbf{x}_t
-\phi^{\mathrm{imp}}(\mathbf{u})
-\sum_{t=1}^T\psi_t(\mathbf{u}_t)\\[2mm]
\mbox{subject to} &
\mathbf{x}_t=\mathbf{x}_{t-1}+\mathbf{u}_t,\quad
\mathbf{x}_0=\mathbf{x}_T=\mathbf{0},\\
&\one^\top\mathbf{x}_t=0,\qquad
\|\mathbf{x}_t\|_1\leq x^{\max},\\
&\mathrm{CVaR}_{\alpha}\big(-\tilde{\mathbf{r}}^\top\mathbf{x}_t\big)\leq r^{\max},
\quad t=1,\ldots,T,
\end{array}
\end{equation}
with variables \(\mathbf{u}_1,\ldots,\mathbf{u}_T\) (dollar trades)
and \(\mathbf{x}_1,\ldots,\mathbf{x}_T\) (dollar holdings).  The objective is forecast return, minus the impact cost
\(\phi^{\mathrm{imp}}\) of Section~\ref{s-model}, minus a
bid/ask transaction cost \(\psi_t(\mathbf{u}_t)=\mathbf{a}^\top|\mathbf{u}_t|\), with
\(\mathbf{a}=(1.0,1.5)\) basis points for the nearby and deferred
legs. The constraint \(\one^\top\mathbf{x}_t=0\) makes
the position a calendar spread with equal dollars long and short, and
\(x^{\max}\) caps its gross notional. The last constraint is a daily risk limit, imposed
for each \(t=1,\ldots,T\): \(\tilde{\mathbf{r}}\) is a one-day return
vector drawn from the empirical distribution of the \(N\) trading
days before the month, and \(\mathrm{CVaR}_{\alpha}\) is the expected
loss over the worst \(1-\alpha\) fraction of those scenarios, so each
day's tail loss is capped at \(r^{\max}\).

We compare three specifications of impact cost:
\begin{enumerate}[leftmargin=1.6em,itemsep=2pt]
\item \emph{baseline}: the separable \(3/2\)-power cost
      \(\sum_{t,i} \ell_i \sigma_{t,i}|u_{t,i}|^{3/2}/v_{t,i}^{1/2}\)
      of \cite{boyd2017multiperiod}, separable across both contracts
      and times, in place of \(\phi^{\mathrm{imp}}\);
\item \emph{self-impact}: \(\mathbf{C}_t=\mathbf{I}\) (no cross-impact)
      and \(\mathbf{G}=\mathbf{I}\) (no transient impact);
\item \emph{transient cross-impact}: the estimated \(\mathbf{C}_t\) and the kernel
      \(\mathbf{G}\).
\end{enumerate}  

\paragraph{Data.}
Prices are the daily settlements of the WTI futures contracts published by the U.S.\ EIA \cite{eia2026futures}. The data covers 2004--2011, when the roll effect was economically
material.  

The return forecast
\(\mathbf{r}_t\) for business day \(t\) of a month is the mean
realized return on business day \(t\) over all prior months of the
sample. 

The dollar volumes \(v_{t,i}\) of \eqref{e-single-impact} are licensed
exchange data and not public, so we simulate them instead.  Let \(v\) denote the total daily dollar volume of both the nearby and deferred contracts in a given month, held constant within the
month.  We simulate \(v\) with public data as
\[
  v
  =\underbrace{q}_{\text{contracts/day}}
  \times\underbrace{1{,}000}_{\text{bbl/contract}}
  \times\underbrace{p_0}_{\text{\$/bbl}},
\]
where \(q\) is the average daily volume of all WTI futures over
the prior calendar year and \(p_0\) is the nearby contract's
settlement price on the month's first trading day.  We read
\(q\) from a CFTC study  \cite{cftc2018tightoil}. The nearby contract's
share of \(v\) is a logistic function of the number of business days
\(e_t\) remaining until its last trading day,
\begin{equation*}
  w_t=\frac{w_0}{1+\exp\big((c-e_t)/s\big)},
\end{equation*}
and the deferred contract takes the rest, i.e. \(v_{t,1}=w_tv\) and
\(v_{t,2}=(1-w_t)v\).  The share starts near \(w_0=0.8\) early in the
month, and falls to half \(c=7\) business days before the last
trading day; with \(s=2\), about half of the \(w_t\)-migration falls inside the five-day roll
window, where 55--62\% of index-roll activity is documented
\cite{irwin2023orderflow}. The Wednesday inventory report also concentrates
liquidity on release days.  We simulate this by scaling \(v\) and \(\sigma_{t,i}\) by 1.08 and
1.02 on Wednesdays, estimated from past daily volume and
return data.

The volatilities and the impact correlation are
estimated by an exponentially weighted moving average (EWMA) of
daily settlement returns, with decay rate \(\lambda=0.97\).  Writing
\(d_m\) for the last trading day of month \(m-1\), we estimate the
contract covariance at \(d_m\) as \(\mathbf{\Sigma}_{d_m}\), and the
volatility of every time period \(t\) in month \(m\) as
\(\sigma_{t,i}=(\mathbf{\Sigma}_{d_m})_{ii}^{1/2}\).  For simplicity, we freeze the estimate
within each month: \eqref{e-roll} is solved once, before the month
begins, so the decision making is based on the information available
at \(d_m\).\footnote{Practitioners would re-solve a receding-horizon
variant of \eqref{e-roll} each day, with daily updated estimates.}
For the risk limit, we set the confidence level \(\alpha=0.95\) and
use each of the \(N=250\) trading days before the month as a
scenario to compute \(\mathrm{CVaR}_{\alpha}\).  We set the
tail-loss cap \(r^{\max}\) to 50 basis points of \(x^{\max}\): the
trader accepts a one-day loss of at most 0.5\% of gross capacity on
the worst 5\% of days.

We set the parameters of temporal decay kernel to \((\tau,\beta)=(1,1/2)\).  The normalizer \(\tau=1\) matches the decision interval, as we consider one decision per day. As we lack available data to learn the
impact scaling \(\ell_i\), we follow the empirics of \cite{boyd2017multiperiod} to
set \(\ell_i=1\), so that trading one day's
volume moves the price by about one day's volatility. 

\paragraph{Metrics.}
We report the average of monthly profit and loss, in basis points of \(x^{\max}=\$1\)B,
\begin{equation}\label{e-pnl}
  \mathrm{net}
  =\underbrace{\sum_{t=1}^{T}\mathbf{r}_t^{\mathrm{realized}}{}^{\top}\mathbf{x}_t}_{\text{alpha}}
  \;-\;\underbrace{\phi^{\mathrm{imp}}(\mathbf{u})}_{\text{impact cost by \eqref{e-causal-cost}}}
  \;-\;\underbrace{\sum_{t=1}^{T}\mathbf{a}^{\top}|\mathbf{u}_t|}_{\text{spread}}.
\end{equation}
We report impact cost in two ways in Table~\ref{t-exp1}, as it cannot be measured directly. Simulated impact cost is calculated by \eqref{e-causal-cost} with the
estimated \(\mathbf{C}_t\) and the kernel \(\mathbf{G}\). 
Predicted impact cost reports each model's
own forecast of \(\phi^{\mathrm{imp}}(\mathbf{u})\) under its cost
specification from the \textbf{Setup} paragraph.  The transient cross-impact model therefore serves two roles: it is
the simulator's cost model and one of the tested models.  
This limitation is standard in impact backtests. 

\begin{figure}[t]
\centering
\includegraphics[width=\textwidth]{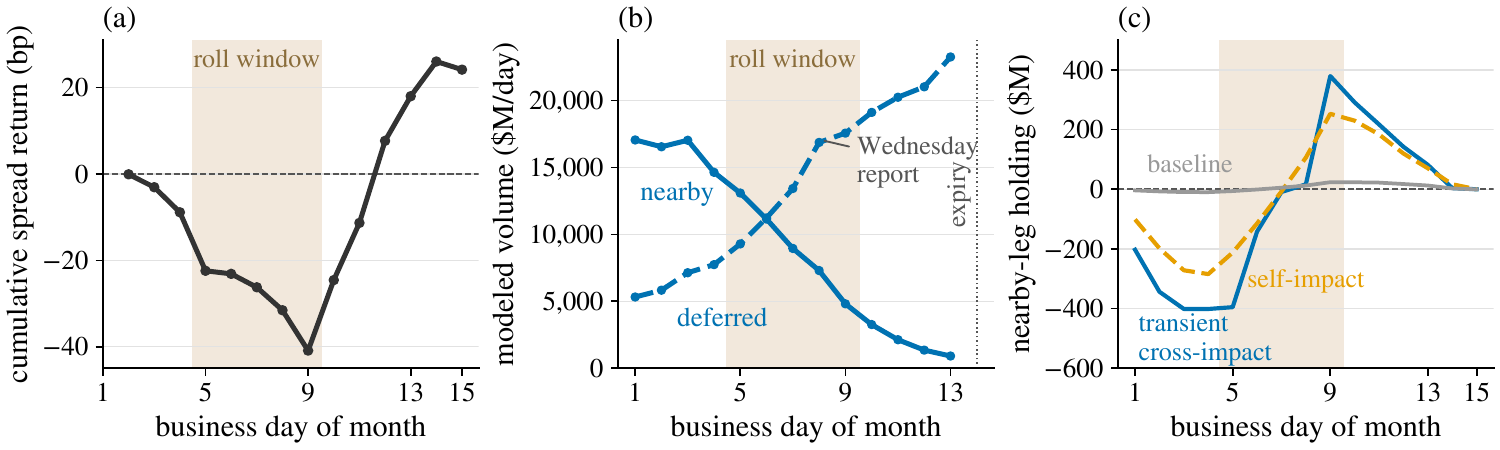}
\caption{Index-roll positioning.  (a)~Averaged cumulative
nearby-minus-deferred spread return by business day of month,
2004--2011; the shaded band is the S\&P GSCI roll window.  (b)~Simulated daily
dollar volumes of the two contracts.  (c)~Averaged nearby contract holding of all three
models.}
\label{f-exp1}
\end{figure}

\paragraph{Results.}
Figure~\ref{f-exp1}(a) shows the price dislocation: averaged over the 96 months, the
nearby-minus-deferred calendar spread falls 9 basis points over business
days two through four, falls another 32 over the roll window, and
rebounds 65 by day 15, consistent with the documented 30--40 basis point
compression and roughly two-week recovery
\cite{irwin2023orderflow}.  Figure~\ref{f-exp1}(b) shows the
liquidity over time: volume migrates from nearby to deferred
contract, with the weekly ripple of the Wednesday inventory report.  Figure~\ref{f-exp1}(c) shows the trading
positions: all short the spread ahead of the roll
window, flip long near its end, and are flat before the nearby's
liquidity dries up.  The models differ in size and timing.

\begin{table}[t]
\centering
\small
\begin{tabular}{lrrrrr}
\toprule
 & \multicolumn{1}{c}{alpha} & \multicolumn{2}{c}{impact cost} &
 \multicolumn{1}{c}{spread} & \multicolumn{1}{c}{net} \\
 \cmidrule(lr){3-4}
policy & & predicted & simulated \eqref{e-causal-cost} & & \\
\midrule
baseline \cite{boyd2017multiperiod} & 0.6 & 1.0 & 0.0 & 0.2 & 0.4 \\
self-impact            & 7.7 & 6.9 & 2.7 & 2.7 & 2.3 \\
transient cross-impact & 12.9 & 5.9 & 5.9 & 3.9 & 3.1 \\
\bottomrule
\end{tabular}
\caption{Index-roll positioning, monthly averages over 2004--2011,
in basis points of the gross limit \(x^{\max}=\$1\)B; the alpha,
impact, spread, and net columns are defined in \eqref{e-pnl}.}
\label{t-exp1}
\end{table}

Figure~\ref{f-exp1}(c) illustrates that the self-impact model and
the baseline \cite{boyd2017multiperiod} overprice their impact.  With estimated correlation \(0.96\), the
price impact of the long and short legs partially cancel, but neither
model captures this netting. The baseline is highly conservative and barely trades, as it charges each contract independently at the \(3/2\) power, so the
spread trade looks expensive and gets no netting.

\section{Conclusion}

We proposed a convex quadratic model of price impact
cost in which a positive semidefinite matrix couples trades across
contracts and a power-law kernel couples trades across dates.  We
showed that the resulting cost admits no price manipulation.  We
empirically demonstrated our model's usefulness, and the overpricing
of existing models, in relative-value trading of crude oil futures.

\section*{Acknowledgements}
We gratefully acknowledge support from the Office of Naval Research under award N000142412306, Air Force Office of Scientific Research under award FA9550-26-1-0012, the Alfred P. Sloan Foundation, the Stanford Institute for Human-Centered Artificial Intelligence (HAI), and from IBM Research as a founding member of Stanford Institute for Human-centered Artificial Intelligence. 

\newpage
\printbibliography

\newpage
\appendix
\section{Model preserves diagonals but not convexity}
\label{a-arith}
The construction \eqref{e-A} couples time periods through the symmetric matrix
$\mathbf{A}_{ts}=G_{ts}\mathbf{\Lambda}_t^{1/2}\mathbf{\Lambda}_s^{1/2}$. Here we show that the true, causal
matrix \(G_{ts}\mathbf{\Lambda}_s\) fails to be positive semidefinite. The counterexample uses $n=1$,
$T=2$. Write $\Lambda_1 = \gamma_1$, $\Lambda_2 = \gamma_2$, and let
$\rho = k(a_2-a_1)\in(0,1]$, so the off-diagonal weight is $G_{12} =
\rho/2$. A model in which the trade at $s$ displaces prices in
proportion to the liquidity prevailing at $s$ gives cross terms
$k(a_t-a_s)\,x_t^T\Lambda_s x_s$ for $t>s$, i.e., 
\[
\mathbf{M} = \begin{bmatrix}
\gamma_1 & \dfrac{\rho\gamma_1}{2} \\[2mm]
\dfrac{\rho\gamma_1}{2} & \gamma_2
\end{bmatrix},
\qquad
\det \mathbf{M} = \gamma_1\gamma_2 - \frac{\rho^2\gamma_1^2}{4},
\]
which is negative whenever $\gamma_2 < \rho^2\gamma_1/4$: an illiquid
period followed closely by a liquid one. 

Also see \cite{fruth2014optimal} for price manipulation in
propagator models with time-varying liquidity.

\end{document}